\documentclass[12pt,A4]{article}
\DeclareSymbolFont{AMSb}{U}{msb}{m}{n}
\DeclareSymbolFontAlphabet{\Bbb}{AMSb}

\usepackage{latexsym}
\usepackage[dvips]{graphicx}
\usepackage{amstext,amsfonts,amsmath,amsthm,graphicx,amssymb,amscd,epsfig}
\usepackage{enumerate}
\usepackage{psfrag} 
\usepackage[all,2cell,ps]{xy}

\input epsf

\newtheorem{teo}{Theorem}[section]

\newtheorem{lem}[teo]{Lemma}
\newtheorem{cor}[teo]{Corollary}
\newtheorem{prop}[teo]{Proposition}
\newtheorem{rem}[teo]{Remark}

\newtheorem{nota}[teo]{Notation}

\newcommand{\norm}[1]{\left\Vert#1\right\Vert}
\newcommand{\abs}[1]{\left\vert#1\right\vert}

\newcommand{\R}{\mathbb{R}}
\newcommand{\Z}{\mathbb{Z}}
\newcommand{\N}{\mathbb{N}}
\newcommand{\no}{\noindent}

\newcommand{\dis}{\displaystyle}

\def\be{\begin{equation}}
\def\ee{\end{equation}}
\def\bq{\begin{eqnarray}}
\def\eq{\end{eqnarray}}
\def\beq{\begin{eqnarray}}
\def\eeq{\end{eqnarray}}
\def\ba{\begin{array}}
\def\ea{\end{array}}
\def\bi{\begin{itemize}}
\def\ei{\end{itemize}}

\newcommand{\spec}{\operatorname{{Spec}}}

\newcommand{\sr}{\operatorname{{sr}}}

\newcommand{\Emb}{\textrm{Emb}}

\title{\bf Planar Embeddings with a Globally Attracting Fixed Point}
\author{Bego\~na Alarc\'on \thanks{Work  partially
supported in part by MECD Grant \# HBE2002-0006, and by MECD Grant
\# MTM2004-03244.}, V\'ictor Gu\'i\~nez
\thanks{Supported in part by FONDECYT Grant \# 1050036.} and
Carlos Gutierrez \thanks{The third author was supported in part by
FAPESP Grant \# TEM\'ATICO 03/03107-9, by CNPq Grant \#
306992/2003-5 and Project CAPES--MECD grant 071/04 and HBP2003--0017
.}}
\date{{}}

\begin{document}
\maketitle


\begin{abstract}
We consider sufficient conditions which guarantee that an embedding
from the plane $\mathbb{R}^2$ into itself has a unique fixed point.
We study sufficient conditions which imply the appearing of a
globally attracting fixed point for such an embedding.
\end{abstract}


\section{Introduction}

\bigskip

This work deals with embeddings from the plane $\mathbb{R}^2$ into
itself. It was motivated by questions about the existence of unique
fixed points, as well as  questions about stability. (See
\cite{Szego}, \cite{Cima-Manosa} \cite{Libro-Jacobian-Conjectur},
\cite{Ortega_degree} and \cite{Robinson}.) It was, in particular,
inspired by the following:

\medskip

\no {\bf DMY Question (Discrete Markus-Yamabe Question)
\cite{Cima-Manosa, Libro-Jacobian-Conjectur}} Let $f:\R^2\rightarrow
\R^2$ be a $C^1-$map such that $f(0)=0$ and $\spec(f)\subset
B(0,1)$. Is $0$ a global attractor for the discrete dynamical system
generated by $f$?

\medskip

To state our results, we shall need the following definitions.

\medskip

Let $f:\mathbb{R}^n \rightarrow
 \mathbb{R}^n$ be a differentiable map. Denote by $Spec(f)$ the set of
eigenvalues of $\; Df_{p} \text{ , for all } p\in \mathbb{R}^n.$

\medskip

 Let $f: \mathbb{R}^n \rightarrow
\mathbb{R}^n$ be a continuous map. We say that $\gamma$ is an
$f-$invariant ray if $\gamma$ is a smooth embedded curve with no self-intersections starting at $0$ and going to
infinity such that $f(\gamma)\subseteq \gamma$. It is assumed that
$0 \in \gamma$.

\medskip

Let $\Emb(\R^2)$ denote the set of topological embeddings $f:\R^2
\rightarrow \R^2$. We will denote by $\Emb_{+}(\R^2)$ the subset of
the orientation preserving maps of $\Emb(\R^2)$ and by
$\Emb^{d}(\R^2)$ the subset of differentiable maps of $\Emb(\R^2).$
Finally, let $\Emb_{+}^{d}(\R^2)=\Emb^{d}(\R^2)\cap\Emb_{+}(\R^2)$

\medskip

Let $f\in \Emb(\R^2)$ and let $p\in \R^2$. We define the
$\omega-$limit set of $p$, $\omega(p)$,  the set (possibly empty)
of points $x\in \R^2$ for which there exists a sequence of
positive integers $n_k \rightarrow \infty$ such that
$f^{n_{k}}(p)\rightarrow x$ as $k \rightarrow \infty$.

\medskip

We will say that $x\in \R^2$ is a
non-wandering point of $f\in \Emb(\R^2)$ if for every neighborhood
$U$ of $x$, there exist an integer $n>2$ and a point $q\in U$ such
that $f^n(q)\in U$. We denote by $\Omega(f)$ the set of
non-wandering points of $f$.

\medskip

%
%
%
Let $f\in \Emb(\R^2)$
\begin{itemize}
\item We say that $0\in\mathbb{R}^2$ is
a local attractor (resp. a local repellor) for $f$ if there exists a
compact disc $D,$ contained in the domain of definition of $f$
(resp. contained in the domain of definition of $f^{-1}$), which is
a neighborhood of $0$ such that $f(D)\subset \mbox{\emph{Int}}(D)$
(resp. $f^{-1}(D)\subset \mbox{\emph{int}}(D)$) and
$\cap_{n=1}^\infty f^n(D) = \{0\}$ (resp. $\cap_{n=1}^\infty
f^{-n}(D) = \{0\}.$ \par

Note that if $f\in \Emb^{1}(\R^2)$ has a hyperbolic attractor at
$0$, our definition coincides with the classical one  \cite{Palis}.

\item  We shall say that $0\in\mathbb{R}^2$ is
a global attractor for $f$ if $0$ is a local attractor and for all
$p \in \mathbb{R}^2$, $\omega(p)=\{0\}$.
\end{itemize}

\medskip

If $f:\mathbb{R}^2\cup\{\infty\} \to \mathbb{R}^2\cup\{\infty\}$
is a homeomorphism of the Riemann Sphere, with $f(\infty)=\infty,$
we may similarly define when $\infty$ is either an atractor or a
repellor.

\medskip

A compact set $W\subset \R^2$ is a window of $f\in \Emb(\R^2)$ if
for all $p\in \R^2$ there exists a $n_{0}\in \N$ such that
$f^{n_{0}}(p)\in W$. We will say that $f\in \Emb(\R^2)$ is a
dissipative map if there exists a window for $f$.

\medskip
Assume that $f\in \Emb(R^2)$ such that $f(p)=p$. The fixed point
$p$ is stable in the sense of Lyapunov (or Lyapunov stable) if
every neighborhood $U$ of $p$ contains another neighborhood $V$
such that, for each $n\geq 0$, $f^n(V)$ is well defined and
$f^n(V)\subset U$.

\medskip

It was proved in \cite{Cima-Manosa} that the DMY Question had a
positive answer for polynomial diffeomorphisms of $\R^2$ but was
false even for rational diffeomorphism of $\R^2$. Therefore, we
wondered if the DMY Question had a positive answer, for smooth
diffeomorphisms of $\R^2$ under the additional assumption
\begin{itemize}
\item[(a)] $\infty$ is a repellor.
\end{itemize}
We show that this assumption is not good enough either: there exists
a smooth diffeomorphism $f:\R^2\rightarrow \R^2$ having an order
four periodic point and such that $\infty$ is a repellor, $f(0)=0$
and $\spec(f)\subset B(0,1)$. (See Theorem~\ref{exemplo}).

\medskip

Since every global diffeomorphism $f$ of $\R^2$ having $0$ as a
global attractor satisfies
\begin{itemize}
\item[(b)] $f$ has an $f-$invariant ray
\end{itemize}
(see Proposition~\ref{infinito-repulsor} where we show that for such
an $f,$ there is a foliation of $\mathbb{R}^2\setminus{0}$ by
$f-$invariant rays), we study the MYD Question under the additional
assumptions (a) and/or (b) obtaining the following five results:

\medskip

\no {\bf Theorem~\ref{alternativa-omega}.} Let $f\in \Emb^d_+(\R^2)$
be such that $ f(0) = 0 $, and that it has an $ f-$invariant ray $
\gamma $. If $\; \spec(f)\cap [1,1+\epsilon)=\emptyset,$ for some
$\epsilon > 0$, then  $\Omega(f) \subset \gamma$ and  either
$\omega(p)=\{0\}$ or $\omega(p)=\emptyset$, for all $p\in \R^2$.

\medskip

\no {\bf Theorem~\ref{cor2curvinvsinptosperiodicos}} Let $f\in
\Emb^d_+(\R^2)$ be such that $ f(0) = 0 $, and that it has an $
f-$invariant ray $ \gamma $. If $\; \spec(f)\subset B(0,1)$, then
$0$ is a local attractor and $\Omega(f)=\{0\}$.

\medskip

\no {\bf Theorem~\ref{chileno}} Let $f:\R^2\rightarrow \R^2$ be a
differentiable map such that $f(0)=0.$ Then:
\begin{enumerate}
\item[(a)] If $\spec(f)\cap ([1,1+\epsilon)\cup [0,\epsilon))=\emptyset$, for
some $\epsilon>0$, then  $ f$ is an
 embedding with $Fix(f)=\{0\}$.
\item[(b)] If $\spec(f)\subset B(0,1)\setminus [0,\epsilon),$ for some
$\epsilon>0$, then not only is $ f$ an
 embedding with $Fix(f)=\{0\}$, and $ 0$ is a local
 attractor  for $ f $.
 \end{enumerate}

\medskip

\begin{figure}[!htb] \label{contraejemplo}
\begin{center}
\begin{minipage}[b]{0.3\linewidth}
\includegraphics[width=3cm]{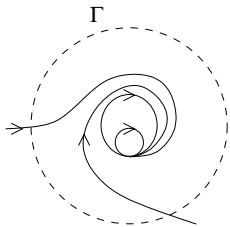}\\
\end{minipage}\hfill
\end{center}
\caption{A vector field which does not have the origin as an
attractor.}
\end{figure}

Relevant to our results, Figure \ref{contraejemplo} shows the phase
portrait of a smooth flow transversal do the unit circle $\partial
D$ such that, the time-one-map $f : \mathbb{R}^2 \to \mathbb{R}^2$,
induced by the flow, which can be assumed to be a smooth
diffeomorphism of $\mathbb{R}^2,$ satisfies  $f(D) \subset
\mbox{Int} \, (D)$ (where $D$ is the unit disc.) and, for all $p \in
\mathbb{R}^2,\;\omega(p)=\{0\}$. Observe that any trajectory of the
flow meeting $\partial D$ is an $f-$invariant ray, while the origin
 $0$ is not a local attractor. For further examples, see \cite{Szego}.

\medskip

\medskip


\no {\bf Theorem~\ref{teoEmbedding}.} Let $ f\in \Emb(\R^2)$ be a
dissipative map such that $ f(0) = 0 $, and that it has an
$f-$invariant ray $\gamma$. Under either of the following conditions
\begin{enumerate}
\item[a)] $ f\in \Emb_+(\R^2)$ and  $ \mbox{Fix}(f) = \{0\} $,
\item[b)] $ \mbox{Fix}(f^2) = \{0\} $,
\end{enumerate}
we have that $\omega(p)=\{0\}$, for all $p\in {\mathbb{R}}^2$.
Moreover, if $0$ is a locally Lyapunov stable fixed point, then $0$
is a global attractor for $f$.

\bigskip

\no {\bf Theorem~\ref{teoDiferencial}.}  Let $f\in \Emb(\R^2)$ be of
class $ C^1 $ such that $f(0)=0$, and that it satisfies the
following three conditions:
\begin{itemize}
\item[(i)] There exist real numbers $R>0$ and $0<\alpha<1$ such
that $||Df_{p}\cdot p|| < \alpha \; ||p||,$  for every $||p||>R.$
\item[(ii)] $Spec(f)\subset B(0,1)$. \item[(iii)] There exists a
$f-$invariant ray $\gamma$.
\end{itemize}

Under either of the following two conditions:
\begin{enumerate}
\item[$(1)$] $f$ is orientation preserving
\item[$(2)$] $\spec(f)\cap \R= \emptyset$,
\end{enumerate}
we have that $0$ is a global attractor for $f$.

\medskip

In obtaining our results, we have profited from the main results of
\cite{aag}, \cite{LCal_teoindice}, \cite{Roland}, \cite{Murthy},
\cite{LRoux_migrations}.

\medskip

Section 2 deals with conditions which imply that a continuous map
from $\mathbb{R}^2$ into itself is an embedding and/or has a unique
fixed point. In Section~3, we study dissipative embeddings having a
global attractor. Some important examples are considered in Section
4. Section 5 presents the results which motivated our assumptions.


\section{Embeddings with only one fixed point}






We shall need the following

\medskip

\begin{teo}[Murthy \cite{Murthy}, 1998] \label{Murthy}
Let $f\in \Emb_{+}(\R^2)$. If $f$ has a periodic point $p$, then
the set $\mbox{Fix}(f)$ of the fixed points of $f$ is nonempty.
\end{teo}

\begin{cor}\label{omega-fix}
Let $f\in \Emb_{+}(\R^2)$. If $\Omega \, (f) \neq \emptyset$ then
$\mbox{\emph{Fix}} (f) \neq \emptyset$.
\end{cor}

\begin{proof}
Suppose by contradiction that $\mbox{Fix}(f)=\emptyset.$ By
Theorem~\ref{Murthy}, $\Omega(f)$ does not contain periodic
points. Let $p\in\Omega(f).$ Then there exist a disc $U,$ an
integer $n>2$ and a point $q\in U$ such that $f^n(q)\in U,$ and
for all, $k=1,2,\dots, n-1,$ $U\cap f^k(U)=\emptyset.$ Let $g:
\R^2 \rightarrow \R^2$ be an orientation-preserving embedding
which is a perturbation of $f$ supported in $U,$ such that
$g(f^n(q))=f(q).$ Then $q$ is periodic point of $g$ of period $n.$
Therefore, by Theorem~\ref{Murthy}, $g$ has a fixed point $q_0.$
Certainly $q_0\notin U$ and so $q_0\in \mbox{Fix}(f).$ This
contradiction proves the result.

\end{proof}

\begin{cor} \label{remMurthy}
If $U\subset \R^2$ is an open set homeomorphic to $\R^2$ and $g: U
\rightarrow U$ is an orientation preserving embedding with
$\Omega(g)\neq \emptyset$, then $g$ has a fixed point.
\end{cor}

\begin{proof}
If $H:U \rightarrow \mathbb{R}^2$ is a homeomorphism then we may
apply Corollary \ref{omega-fix} to $f=H \circ g \circ H^{-1}$.
\end{proof}

\medskip

We remark that the version of Theorem \ref{Murthy} for
chain-recurrent-points can be found in \cite{Murthy}.

\medskip

\begin{teo}[Fernandes, Gutierrez and Rabanal \cite{Roland}, 2004] \label{teoroland}
Let $f:\R^2\rightarrow \R^2$ be a differentiable (not necessarily
$C^1$) map. If for some $\epsilon
>0,$ $\; \spec(f)\cap [0,\epsilon)=\emptyset$, then $f$ is an
embedding.
\end{teo}

\begin{cor} \label{propptofijo}
Let $f:\R^2 \rightarrow \R^2$ be a differentiable map such that $
\spec(f)\cap [1,1+\epsilon)=\emptyset, \text{ for some }
\epsilon>0.$ Then \emph{Fix}$(f) $ is either empty or a
one--point--set.
\end{cor}

\begin{proof}
Let $g=f-Id$, where $Id:\mathbb{R}^2 \rightarrow \mathbb{R}^2$ is
the identity map. Certainly  $\lambda \in \spec(f)$ if and only if
$\lambda-1\in\spec(g)$. So there exists an $\epsilon>0$ such that
$\spec(g)\cap[0,\epsilon)=\emptyset$. By Theorem \ref{teoroland},
$\,g$ is injective; this implies the conclusion.
\end{proof}

The following result gives us a condition for a differentiable map be an
embedding having a unique fixed point:

\begin{cor} \label{propptofijo1}
Let $f:\R^2 \rightarrow \R^2$ be a differentiable map such that $f(0)=0
\text{ and } \spec(f)\cap (
[0,\epsilon)\cup[1,1+\epsilon))=\emptyset, \text{ for some }
\epsilon>0.$ Then $f\in \Emb^d(\R^2)$ such that $\emph{Fix}(f) =
\{0\}$.
\end{cor}

\begin{proof}
It follows immediately from Theorem~\ref{teoroland} and
Corollary~\ref{propptofijo}.
\end{proof}

\begin{lem} \label{propcurvinvsinptosperiodicos} Let $f\in
\Emb(\R^2)$ with $ f(0) = 0 $ and having an $ f-$invariant ray $ \gamma $.
Under either of the following conditions
\begin{enumerate}
\item[a)] $f\in \Emb_+(\R^2)$ and $ \mbox{Fix}(f) = \{0\} $,
\item[b)] $ \mbox{Fix}(f^2) = \{0\} $,
\end{enumerate}
we have that $\Omega(f)
\subset \gamma$ and  either $\omega(p)=\{0\}$
or $\omega(p)=\emptyset$, for all $p\in \R^2$.
\end{lem}

\begin{proof}
Suppose a) is satisfied. If $0$ is the only fixed point of $f$ and
there exists a $f-$invariant ray $\gamma$, then $f:\R^2\setminus
\gamma \rightarrow \R^2\setminus \gamma \;$ is a fixed-point-free
orientation preserving embedding.\par

\medskip

As $\R^2\setminus \gamma \subset \R^2$ is an open set homeomorphic
to $\R^2$, by Corollary~\ref{omega-fix},
$\Omega(f_{|_{\mathbb{R}^2\setminus \gamma}})=\emptyset$. Then
$\Omega(f) \subset \gamma$ and since $f|_{\gamma}: \gamma
\rightarrow \gamma$ is a one dimensional homeomorphism with only one
fixed point, then given $p\in \R^2$ either $\omega(p)=\{0\}$ or
$\omega(p)=\emptyset$.

\medskip

If b) is satisfied we may apply a) to $ f^2, $ which is always
orientation preserving, and obtain the required conclusion.
\end{proof}


\begin{teo}\label{alternativa-omega}
Let $f\in \Emb^d_+(\R^2)$ be such that $ f(0) = 0 $, and that it has
an $ f-$invariant ray $ \gamma $. If $\; \spec(f)\cap
[1,1+\epsilon)=\emptyset,$ for some $\epsilon > 0$, then  $\Omega(f)
\subset \gamma$ and  either $\omega(p)=\{0\}$ or
$\omega(p)=\emptyset$, for all $p\in \R^2$.
\end{teo}
\begin{proof}
It follows directly from Lemma \ref{propcurvinvsinptosperiodicos}
and Corollary \ref{propptofijo}.
\end{proof}

\begin{teo} \label{cor2curvinvsinptosperiodicos}
Let $f\in \Emb^d_+(\R^2)$ be such that $ f(0) = 0 $, and that it has
an $ f-$invariant ray $ \gamma $. If $\; \spec(f)\subset B(0,1)$,
then $0$ is a local attractor and $\Omega(f)=\{0\}$.
\end{teo}
\begin{proof}
It follows from the Hartman and Grobman Theorem \cite{Palis} that
$0$ is a local attractor. By applying
Theorem~\ref{alternativa-omega}, we obtain this result.
\end{proof}

\medskip

\begin{rem} Under the condition of Theorem~\ref{teoEmbedding}, if $0$ is stable in
the sense of Lyapunov, we have $ \emph{index}(f,0)=1 $
(\cite{Ortega_degree}).
\end{rem}

\medskip

\begin{teo}\label{chileno} Let $f:\R^2\rightarrow \R^2$ be a
differentiable map such that $f(0)=0.$ Then:
\begin{enumerate}
\item[(a)] If $\, \spec(f)\cap ([1,1+\epsilon)\cup [0,\epsilon))=\emptyset$, for
some $\epsilon>0$, then  $ f$ is an
 embedding with $Fix(f)=\{0\}$.
\item[(b)] If $\, \spec(f)\subset B(0,1)\setminus [0,\epsilon),$ for some
$\epsilon>0$, then not only is $ f$ an
 embedding with $Fix(f)=\{0\}$, and $ 0$ is a local
 attractor  for $ f $.
 \end{enumerate}
\end{teo}
\begin{proof} Item (a) follows from
Corollary~\ref{propptofijo}. Item (b) follows from (a) and from
Hartman and Grobman Theorem (\cite{Palis}).
\end{proof}

The proof  of theorem below can be found in \cite{aag}.

\begin{teo}\label{chileno2} Let $f:\R^2\rightarrow \R^2$ be a differentiable map with
$f(0)=0$. If $\spec(f)\cap \mathbb{R} = \emptyset$, then $f$ is an
embedding with $Fix(f^2)=\{0\}$.
\end{teo}





\medskip

\section{Dissipative embeddings and attractors}

\medskip

We shall need the following result (see also \cite{Martins}):

\begin{prop}(Richeson and Wiseman \cite{RiWi}, 2002). \label{Martins}
Let $X$ be a metric locally compact set and $f:X \rightarrow f(X)$
be an homeomorphism. The following are equivalent:
\begin{enumerate}
\item[(i)] $f$ is a dissipative map. \item[(ii)] For all compact
set $D\subset X$ there exists a window $W$ for $f$ such that
$D\subset int(W)$, $f(W) \subset W$ and $\forall p\in X \; \exists
n_{0}$ such that $f^{n_{0}}(p) \in int(W)$. \item[(iii)] For all
compact set $D\subset X$ there exists a window $W$ for $f$ such that
$D\subset int(W)$ and $f(W) \subset int(W)$.

\end{enumerate}

\end{prop}

\begin{rem}\label{window}
It follows from Proposition right above  that if $f\in \Emb(\R^2)$
is a dissipative embedding,
then there exists a compact set $D$ such that for all $p\in \R^2$
the set $\omega(p)$ is non empty and $\omega(p)\subset D$. Note that
if $f(0)=0$, then $0\in D$.
\end{rem}

\medskip

\begin{teo} \label{teoEmbedding} Let $ f\in \Emb(\R^2)$ be a
dissipative map such that $ f(0) = 0 $, and that it has an
$f-$invariant ray $\gamma$. Under either of the following conditions
\begin{enumerate}
\item[a)] $ f\in \Emb_+(\R^2)$ and  $ \mbox{Fix}(f) = \{0\} $,
\item[b)] $ \mbox{Fix}(f^2) = \{0\} $,
\end{enumerate}
we have that $\omega(p)=\{0\}$, for all $p\in {\mathbb{R}}^2$.
Moreover, if $0$ is a locally Lyapunov stable fixed point, then $0$
is a global attractor for $f$.
\end{teo}

\begin{proof} As $f$ is a dissipative map, there exists a compact set $D\subset \R^2$
such that for all $p\in\R^2$, $\omega(p)\neq \emptyset$ and
$\omega(p)\subset D$. Using Lemma~\ref{propcurvinvsinptosperiodicos}
we conclude that for all $p\in \R^2$, $\omega(p)=\{0\}$. Therefore,
if $0$ is locally stable in the sense of Lyapunov, we may obtain
that $0$ is a global attractor.
\end{proof}

\medskip

\begin{rem} Under the condition of Theorem~\ref{teoEmbedding}, if $0$ is stable in
the sense of Lyapunov, we have $ \emph{index}(f,0)=1 $
(\cite{Ortega_degree}).
\end{rem}

\begin{lem}\label{propbola} Let $f:\R^n \rightarrow \R^n$
be a $C^1-$map such that $f(0)=0$. Suppose that there exist $R>0$
and $0<\alpha<1$ such that
$$\forall \; \norm{p}>R, \; \quad \norm{Df_{p}\cdot p} < \alpha \norm{p}.$$
There exist $S_{0}>R$ and $\alpha < \mu < 1$ such that if $p \in
\mathbb{R}^n$ and $||p|| \geq S_0$ then $||f(p)|| \leq \mu \,
||p||$.
\end{lem}

\begin{proof}

As $\bar{B}(0,R)$ is a compact set, there exists a real number $M>1$
such that $\norm{Df_{q}}\leq M$ for all $q\in \bar{B}(0,R)$.
\par

Let $\mu=\frac{\alpha+1}{2}$ and $S_0=\frac{2(MR-\alpha
R)}{1-\alpha}$. Notice that $0<\alpha<\mu<1$ and $S_0>R$. Given $p
\in \mathbb{R}^n$ with $||p||\ge S_0$, let
$\gamma:[0,||p||]\;\rightarrow\; \mathbb{R}^n$ be defined by
$\gamma(t)=(t\, p)/||p||$. Then
\begin{align*}
||f(p)|| & \leq  \int^R_0 \,||(f \circ \gamma)'\,(t) ||\, dt +
\int^{||p||}_R\, ||(f \circ \gamma)'\,(t)||\, dt\\
& \leq MR+ \alpha\,(||p||-R)\\
& = \displaystyle{\frac{1-\alpha}{2} \,\cdot\,
\frac{2(MR-\alpha R)}{1-\alpha}}+\alpha \, ||p||\\
& \le \displaystyle{\frac{1-\alpha}{2}\, ||p|| + \alpha \,||p|| =
\mu\, ||p||}.
\end{align*}
\end{proof}

\medskip



\begin{teo} \label{teoDiferencial} Let $f\in \Emb(\R^2)$ be of
class $ C^1 $ such that $f(0)=0$, and that it satisfies the
following three conditions:
\begin{itemize}
\item[(i)] There exist real numbers $R>0$ and $0<\alpha<1$ such
that $||Df_{p}\cdot p|| < \alpha \; ||p||,$  for every $||p||>R.$
\item[(ii)] $Spec(f)\subset B(0,1)$. \item[(iii)] There exists a
$f-$invariant ray $\gamma$.
\end{itemize}

Under either of the following two conditions:
\begin{enumerate}
\item[$(1)$] $f$ is orientation preserving
\item[$(2)$] $\spec(f)\cap \R= \emptyset$,
\end{enumerate}
we have that $0$ is a global attractor for $f$.
\end{teo}

\begin{proof}
Suppose that (1) is satisfied.  By item $(i)$ and Lemma
\ref{propbola}, $f$ is dissipative.
Therefore, by items (ii), (iii) and Hartman and Grobman Theorem
\cite{Palis}, we may apply Theorem
\ref{cor2curvinvsinptosperiodicos} to obtain the requested
conclusion: $0$ is a global attractor for $f.$

By Theorem \ref{chileno2}, the proof of this theorem under
assumption (2) is the same as that under the assumption (1); notice
that $f^2$ is always orientation preserving.
\end{proof}


\section{Examples}

The following two examples can be found in \cite{Cima-Manosa}.
\begin{teo}[Szlenk's Example] \label{szlenk}
Let $F:\R^2\rightarrow \R^2$ be defined by
$$F(x,y)=(-\frac{ky^3}{1+x^2+y^2},\frac{kx^3}{1+x^2+y^2}),
\quad \text{ where } \; k\in (1,\frac{2}{\sqrt{3}}).$$ The map $F$
satisfies the following properties:

\begin{enumerate}
\item Set $p=(x,y)\in \R^2$ and let $\lambda$ be an eigenvalue of
$JF(p)$. If $xy=0$ then $\lambda=0$. Otherwise $\lambda \notin \R$
and $\abs{\lambda}<\sqrt{3}k/2$. \item $F^4(
(k-1)^{-1/2},0)=((k-1)^{-1/2},0).$ \item $F$ is injective.
\end{enumerate}
\end{teo}

\begin{teo}[Szlenk--Cima--Gasull--Ma\~nosas's Example]\label{teo-szlenk-cat}
Let $F:\R^2\to\R^2$ be as in theorem~\ref{szlenk} and let $0<a<1.$
Let $G_a:\R^2\to \R^2$ be defined by
$$G_a(x,y) = F(x,y) - a(x,y).$$
If $a$ is small enough, the map $G_a$ is a global diffeomorphism of
$\R^2$ onto itself which satisfies the following properties:
\begin{enumerate}
\item[(a)] for all $x\in \R^2,$ $\spec(G_a)(x)\subset B(0,1)$
\item[(b)] $G_a(0)=0$ and there exists $p\in\R^2\setminus \{0\}$
such that $(G_a)^4(p)=p.$
\end{enumerate}
\end{teo}

\bigskip

\begin{rem} For the Szlenk's map $ F,$ there are points $ p $ such that
$ \lim_{n \to \infty} \; F^n(p) \; = \; \infty $ and therefore $
\infty $ is not a repellor. Also there is numerical evidence that
for the Szlenk--Cima--Gasull--Ma\~nosas's map $ G_a,$ $ \infty $ is
not a repellor.
\end{rem}

\bigskip

The next example shows that DMY Question has a negative answer, for
smooth diffeomorphisms of $\R^2$, under the additional condition
that $ \infty $ is a repellor.

\begin{teo}\label{exemplo}
There exists a smooth diffeomorphism $f:\R^2\rightarrow \R^2$
having an order four periodic point and such that $\infty$ is a
repellor, $f(0)=0$ and $\spec(f)\subset B(0,1).$
\end{teo}

\no We introduce some notations and Lemmas that we will use in the proof of
Theorem \ref{exemplo}.

\begin{nota} If $A:\R^2 \rightarrow \R^2$ is a linear map, we denote by
$\sr\,(A)$ the espectral radius of $A$:
$$\sr\,(A)=\max \, \{| \lambda |:\lambda \;\mbox{is an eigenvalue of} \;
A\}$$ As usual, the norm $||A||$ of $A$ is
$$||A||=\sup \, \{||Av||:||v||=1\}.$$
It is well known that
$$\sr \,(A) \,\leq\, ||A||.$$
If $F:\R^2 \rightarrow \R^2$ is a diffeomorphism we denote by
$$\sr\, (F)= \; \sup\, \{ \sr\, (DF_{(p)}): p \in \R^2\}.$$
\end{nota}

\begin{lem}\label{lem-prev1}
Given $ R > 0 $, $ C > 0 $ and $ 0 < \varepsilon <  \dis{\frac{1}{8C}} $,
there exists a smooth function
$$ \phi:[0,\infty) \,\rightarrow\, \left[\dis{\frac{1}{2C}},1\right] $$
which verifies
\begin{enumerate}
\item[1)] $ \phi (r)=1 $, for all $ r \in [0,R] $,
\item[2)] $ \phi'(r) \leq 0 $ and $ |\phi'(r) \cdot r|\, < \varepsilon $
for all $ r \in [0,\infty) $,
\item[3)] there is an integer $ N $ such that $ \phi (r) = \dis{\frac{1}{2C}} $,
for all $ r > N $.
\end{enumerate}
\end{lem}
\begin{proof}
First we  define $ \phi (r) = 1 $, for all $r \in [0,R] $. Hence, we consider
the sequence $ (S(n)) $ defined for all positive integer $ n $ by

  $$ S(n) \; = \; 1 - \frac{\varepsilon}{8}\left(\frac{1}{R+1} + \frac{1}{R+2} + \cdots + \frac{1}{R+n}\right) \, .$$
 Observe that $ S(n) \to -\infty $ as $ n \to \infty $. Let $ N $ be the
integer such that $ S(n) > \dis{\frac{1}{2C}}
 $ for $ n = 1, \cdots, N - 1 $, and $ S(N) \leq  \dis{\frac{1}{2C}} \, \cdot
$

 \no For all positive integer $n \geq 1 $ we define

  $$\left\{\begin{array}{lcr}
  \phi(R+n)=S(n), & \mbox{if} & 1 \leq n \leq N-1 \; , \\
  \quad & \quad & \quad\\
  \phi(R+n)=\dis{\frac{1}{2C}}, & \mbox{if} &  n \geq N \;
\cdot
  \end{array}\right.$$

\no For each integer $ n $ such that $ 1 \leq n \leq N - 1$ we
consider a smooth map $ \phi \, : \, [R+n , R+n+1] \rightarrow [S(n+1),S(n)] $
which is flat at $ R + n $ and $ R + n +1 $ and defined as above and such that
$$ 0 \leq - \phi'(r) = |\phi'(r)| \leq \frac{2(S(n) - S(n+1))}{(R + n + 1) -
(R + n)} = \dis{\frac{\varepsilon}{4(R+n+1)}} \, \cdot $$

\no Finally, for $ r > R + N $ we define $ \phi(r) = \dis{\frac{1}{2C}}
\, \cdot $

\no Then the map $ \phi $ verifies the conditions because, for $ n = 1, \cdots , N
$ and
$ R + n - 1 \leq r \leq R + n $ we have
$$   |\phi'(r) \cdot r| \leq  \dis{\frac{\varepsilon}{4(R+n)}} \cdot (R + n)
\; = \; \dis{\frac{\varepsilon}{4}} \, < \,\varepsilon \, . $$

\end{proof}
\begin{lem}\label{lem-prev2}
Let $M_2(\R)$ be the space of $2 \times 2$ real
matrices and let $ C $ be a positive constant. Denote by $ \mathcal{A} $ the
compact set consisting of $ A \in M_2(\R) $
such that $ ||A|| \leq C $ and $ \sr(A) \leq 0 \cdot 9 $.
There exists $ \dis{\frac{1}{8C}} > \varepsilon > 0 $ such that if
$$ B=\left(\begin{array}{cc}
b +\varepsilon_1 & \varepsilon_2\\
\varepsilon_3 & b +\varepsilon_4
\end{array}\right) $$
satisfies: $\dis{\frac{1}{2C}} \leq b \leq1$, and $\max
\,\{\,|\varepsilon_1|, |\varepsilon_2|, |\varepsilon_3|,
|\varepsilon_4|\}< \varepsilon$,
then
$$ \sr \, (BA) \leq 0 \cdot 95 \, ,$$
for all $ A \in \mathcal{A} $.
\end{lem}

\begin{proof} Writting $ B = b \, I \, + \, E $, we have $ BA = b A + EA $ and
$$ \det(BA - \lambda I) \, = \, b^2 \det\left(A  + \frac{1}{b} \, EA -
\frac{\lambda}{b} \, I \right) \, , $$
which implies
$$ \sr(BA) \; = \; a \, \sr\left(A  + \frac{1}{b}\, EA \right) \, . $$
Then, given $ \delta > 0 $ there exist $ 0 < \varepsilon < \frac{1}{8C} $ such
that
$$ \sr\left(A  + \frac{1}{b}\, EA \right)  \; \leq  \; \sr(A) \, + \,
\frac{\delta}{b} \; , $$
that is, such that
$$ \sr(BA) \; \leq \; b \, \sr(A) \, + \, \delta \, . $$
Hence, it sufficient to consider $ 0 < \delta < 0,95 \, - \, b \cdot 0,9 \,
$.
\end{proof}

\begin{proof}[\bf{Proof of Theorem \ref{exemplo}}]
Consider the Szlenk--Cima--Gasull--Ma\~nosas's Example
$$ G_a(x,y) \; = \;
(-\frac{ky^3}{1+x^2+y^2},\frac{kx^3}{1+x^2+y^2}) \; - \; a \, (x,y) \, , $$
with $ k\in (1,\frac{2}{\sqrt{3}})  $ and $ 0 < a < a_0 < 1 $ such that
$ G_a $ verifies Theorem \ref{teo-szlenk-cat}. Also consider  a constant
$ C > 0 $ such
that for all  $ (x,y) \in \R^2 $,
\begin{equation}\label{eq1}
||D\, G_a(x,y)|| \leq C \, .
\end{equation}
Since, $ \sr\, (G_0) < \sqrt{3}k/2 $ and
$ \spec(G_a)(x,y) = \spec(G_0)(x,y) - a $, for all  $ (x,y) \in \R^2 $, there
exist $ 0 < a_0 < 1 $
and $ 1 < k_0 < \frac{2}{\sqrt{3}} \cdot (0,88) \, (\approx 1,01614) $, such
that
\begin{equation}\label{eq2}
\sr \, (G_a) \leq 0,9 \, , \quad \textrm{for all} \quad  0 < a < a_0
\quad \textrm{and} \quad 1 < k < k_0 \, .
\end{equation}
Asociated to the numbers $ R= \dis{\frac{2}{\sqrt{k-1}}} $, $\, C > 0 \,$ as in
(\ref{eq1}) and $ 0 < \varepsilon < \dis{\frac{1}{8C}} $ as in Lemma
\ref{lem-prev2}, consider the
smooth function
$$ \phi:[0,\infty) \,\rightarrow\, \left[\dis{\frac{1}{2C}},1\right] $$
of Lemma \ref{lem-prev1}.
We will prove that
$$ f = H \circ G_a  :\R^2 \rightarrow \R^2 \, , $$
with $ a $ (and $ k $) as in (\ref{eq2}), and with $ H $ defined by
$$ H(x,y)=(\phi(\sqrt{x^2+y^2})x,\phi(\sqrt{x^2+y^2})y) \, ,$$
verifies our Theorem.

\no 1) To prove $ f $ is a smooth diffeomorphism we show that $ H $ is a
smooth diffeomorphism. If $ r=\sqrt{x^2+y^2} $,
we have
$$ JH(x,y) \; = \; \left(\begin{array}{ccc}
\phi(r)+\phi'(r) \cdot \dis{\frac{x^2}{r}} & \quad & \phi'(r) \cdot
\dis{\frac{xy}{r}}\\
\quad & \quad & \quad\\
\phi'(r) \cdot \dis{\frac{xy}{r}} & \quad & \phi(r)+\phi'(r) \cdot
\dis{\frac{y^2}{r}}
\end{array}\right)$$

\no Since
$$ \phi(r) > \frac{1}{2C} > 4 \varepsilon > \varepsilon > |\phi'(r) \cdot r|
\, , $$
implies
$$ \det \, DH(x,y) \,=\, \phi(r)\,(\phi(r)+\phi'(r) \cdot r) \,=\,
\phi(r)\,(\phi(r) -  |\phi'(r) \cdot r)|) > 0 \, , $$
and
$$ \lim_{|(x,y)| \to \infty} \; |H(x,y)| \; = \; \infty \, , $$
the map $ H $ is a global diffeomorphism.

\no 2) To prove $ \sr(f) < 1 $ first observe that $ D G_a (x,y) $ is in the
set $ \mathcal{A} $ of Lemma \ref{lem-prev1}, for all $ (x,y) $. Also, if
$ r = |Ga(x,y)| $ and $ (u,v) = Ga(x,y) $ we have
$$ \dis{\frac{1}{2C}} \leq \phi(r) \leq 1 \, , $$
and
$$ \max\{|\phi'(r) \cdot \dis{\frac{u^2}{r}}|,|\phi'(r) \cdot
\dis{\frac{uv}{r}}|,|\phi'(r) \cdot
\dis{\frac{v^2}{r}}| \; \leq \;  |\phi'(r) \cdot r| < \varepsilon \, .
$$
Then, by Lemma \ref{lem-prev2}, we obtain
$$ \sr(f) \; \leq \; 0.95 \, . $$

\no 3) Clearly $ f(0,0) = (0,0) $. Moreover, since $ p = (R/2,0) $
is a hyperbolic order four periodic point of $ G_0 $, then not only
the diffeomorphism $ G_a $ has a  hyperbolic order four periodic
point in the puncture ball $ B(0,R)-\{0\} $ but also $ f $.

\no 4) Finally
$$ |f(x,y)| \, = \, \frac{|G_a(x,y)|}{2C} \; \leq  \; \frac{1}{2} \, |(x,y)|
\, , $$
for all $ |(x,y)| \geq R + N $ with $ N $ as in item 3)  of Lemma
\ref{lem-prev1}. Therefore, $ \infty $ is a repellor for $ f $ and
the proof is complete.
\end{proof}

\section{Final remarks}

The following result might have motivated the DMY Question:
\begin{prop} Let $f:\R^{n}\rightarrow \R^{n}$ be a $C^1-$map such
that $f(0)=0$. If for all $p\in \R^n,$ $\;\norm{Df_{p}}< 1$, then
$0$ is a global attractor for $f$.
\end{prop}
\begin{proof}
Let $q\in \R^n$ and $S=\norm{q}$. As $f$ is a $C^1-$map, and
$\bar{B} (0,S)$ is a compact set there exists a real number $M>0$
such that for all $p\in \bar{B}(0,S),\;$ $\norm{Df_{p}}\leq M<1$.
Hence $f(\bar{B}(0,S))\subseteq B(0,S)$ and so $f|_{\bar{B}(0,S)}$
is a $M-$contraction. Therefore $f^{k}(q)\rightarrow 0$ as
$k\rightarrow \infty$. As $q \in \mathbb{R}^n$ is an arbitrary
point, $0$ is a global attractor for the discrete dynamical system
generated by $f$.
\end{proof}
Next proposition justifies one of our main assumptions.

\begin{prop}\label{infinito-repulsor}
Let $ f: (\R^2,0) \rightarrow (\R^2,0) $ be a $C^r-$diffeomorphism,
with $r\ge 1,$ having $ 0 $ as a global attractor. Then there exists
a $C^r-$foliation of $\, \R^2-\{0\} $ by $ f-$invariant rays.
\end{prop}

\begin{proof} Given $ 0 < r_1 < r_2 $ and $ r > 0 $ let denote
$ A(r_1, r_2) = \{p \in \R^2 : r_1 \leq \norm{p} \leq r_2  \} $, $
B(r) = \{p \in \R^2 : \norm{p} < r \} $ and $ S(r) = \{p \in \R^2 :
\norm{p} = r \} $. Also, given a simple closed curve $ \gamma
\subset \R^2 $, we denote by $ \gamma^- $ (resp. $ \gamma^+ $) the
union of $ \gamma $ with the bounded (resp. unbounded) connected
component of $ \R^2 -\gamma.$

 Let $ D $ be a smooth compact disk such that $ 0 \in \mbox{int}(D),$
 $ f(D) \subset \mbox{int}(D) $ and $ \bigcup_{n\in \Z} f^{n}(D)=
\R^2 $. Denote by $ C = \partial D $ the boundary of $ D $.
 Let $ V(C) $ be a small smooth tubular neighborhood of $ C $ such
that $ f(V(C)) \cap V(C) = \emptyset $.  For $ \varepsilon
> 0 $ small, consider a smooth embedding
$$
g_1: A(2-2\varepsilon,2 + 2\varepsilon)\to V(C) \quad
\mbox{with}\quad g_1(S(2)) = C.
$$
and define the smooth embedding $ g_2 = f \circ g_1 \circ T:
A(1-\varepsilon,1 + \varepsilon)\to(V(C)) $, where $ T(x,y) = 2(x,y)
$. In this way,

\bigskip

\noindent (1) \, restricted to $A(2-2\epsilon, 2+2\epsilon),$ the
map $f\circ g_1$ equals to  $g_2\circ T^{-1}.$

\bigskip

We claim that

\bigskip

\noindent (2) \, there exists a smooth embedding $H_0:
A(1-\varepsilon,2 + 2\varepsilon)\to \mathbb{R}^2$ which is an
extension of both $ g_1 $ and $ g_2 $.

\bigskip

In fact, let $ G_1 : S(2- 2\varepsilon)^+ \rightarrow g_1(S(2-
2\varepsilon))^+ $ and $ G_2 : S(1 + \varepsilon)^- \rightarrow
g_2(S(1 + \varepsilon))^- $ be smooth diffeomorphisms which are
extensions of $ g_1 $ and $ g_2 $, respectively (see Theorems 8.3.3
and and 8.1.9 of \cite{Hirsch}). Using Theorem 8.3.2 of
\cite{Hirsch} we can find a smooth diffeomorphism $ H :
\R^2_{\infty} \rightarrow \R^2_{\infty} $ diffeotopic to the
identity of $\mathbb{R}^2 $ which restricted to $S(2-
2\varepsilon)^+ $ and $ S(1 + \varepsilon)^-$ coincides with $G_1$
and $G_2,$ respectively. The proof of claim (2) is obtained by
defining $H_0$ to be equal to the restriction of $H$ to
$A(1-\varepsilon,2 + 2\varepsilon).$

\bigskip

\noindent On the annulus $ A(1-\epsilon, 2+\epsilon) $ we consider
the foliation $\mathcal{L}_0$ whose leaves are  the linear segments
$ r_p = \{t p : 1-\epsilon \leq t \leq 2+\epsilon \} $ with $ p \in
S(1) $.  We claim that

\bigskip

\noindent (3) \, If $\ell \subset A(2-2\epsilon, 2+2\epsilon),$ is
an arc of a  leaf of $\mathcal{L}_0,$ then  both $H(\ell)$ and
$f(H(\ell))$ are arcs of the same leaf of the foliation
$H(\mathcal{L}_0).$

\bigskip

\noindent In fact, using (1), we obtain that
$$ f\circ H (\ell) = f\circ g_1 (\ell) = g_2\circ T^{-1}(\ell) =
H(T^{-1}(\ell).$$ The claim follows from this and from the fact that
$\ell$ and  $T^{-1}(\ell)$ are arcs of the same leaf of
$\mathcal{L}_0.$

\bigskip

\noindent  Then, by (2) and(3), $ \mathcal{F} = \bigcup_{k\in \Z}
f^{k}(H(\mathcal{L}_0)) $ is a smooth foliation on $\R^2-\{0\}$ by
$f-$invariant rays.

\end{proof}




\end{document}